\documentclass[conference]{IEEEtran}
\IEEEoverridecommandlockouts
\usepackage{cite}
\usepackage{amsmath,amssymb,amsfonts}
\usepackage{algorithmic}
\usepackage{graphicx}
\usepackage{textcomp}
\usepackage{xcolor}
\def\BibTeX{{\rm B\kern-.05em{\sc i\kern-.025em b}\kern-.08em
    T\kern-.1667em\lower.7ex\hbox{E}\kern-.125emX}}
    \newtheorem{thm}{Theorem}
\newtheorem{lem}[thm]{Lemma}
\newtheorem{rem}[thm]{Remark}

\newtheorem{pf}{Proof}
\IEEEoverridecommandlockouts
\IEEEpubid{\makebox[\columnwidth]{ 978-1-6654-8724-5/22/\$31.00~\copyright{}2022 IEEE \hfill} \hspace{\columnsep} \makebox [\columnwidth]{ }}

\begin{document}

\title{Bilinear optimal stabilization of a non-homogeneous Fokker-Planck equation}

\author{\IEEEauthorblockN{ K. Ammari}
\IEEEauthorblockA{\textit{dUR Analysis and Control of PDE's} \\
\textit{Faculty of Sciences of Monastir, University of Monastir}\\
 5019 Monastir, Tunisia \\
kais.ammari@fsm.rnu.tn}
\and
\IEEEauthorblockN{ M. Ouzahra}
\IEEEauthorblockA{\textit{Laboratory MMPA, Department of mathematics $\&$ informatics,} \\
\textit{ ENS. University of Sidi Mohamed Ben Abdellah}\\
 Fes, Morocco \\
mohamed.ouzahra@usmba.ac.ma}
\and
\IEEEauthorblockN{ S. Yahyaoui}
\IEEEauthorblockA{\textit{Laboratory MMPA, Department of mathematics $\&$ informatics,} \\
\textit{ENS. University of Sidi Mohamed Ben Abdellah}\\
Fes, Morocco \\
soufiane.yahyaoui@usmba.ac.ma}
}

\maketitle

\begin{abstract}
In this work, we study the bilinear optimal stabilization of a non-homogeneous Fokker-Planck equation. We first study the problem of optimal control in a finite-time interval and then focus on the case of the infinite time horizon. We further show that the obtained optimal control guarantees the strong stability of the system at hand. An illustrating numerical example is given. 
\end{abstract}


\section{Introduction}
In this paper we consider the following non-homogeneous bilinear system:
{\footnotesize
\begin{equation}\label{P3-sys obj}
\begin{cases}
\dot{y}(t,x)= \Delta y(t,x) +\sum_{i=1}^{N} u_i (t) \left(\frac{\partial y}{\partial x_i} (t,x)+ b_i \right),  (t,x) \in (0,T) \times  \Omega  \\
y(t,x)=0, \, \ \ \ \ \ \ \ \ \ \ \ \ \ \ \ \ \ \  \  \ \ \ \ \ \ \ \ \ \ \ \ \ \ \ \ \ \ \ \ \ \ \ \ \ \ \ \ (t,x) \in (0,T) \times \partial\Omega\\
    y(0,x)=y_0  \in L^2(\Omega),  \ \ \ \ \ \ \ \ \ \ \ \ \ \ \ \ \ \  \ \ \ \ \ \ \ \ \ \ \ \ \ \ \ \ \ \ \ \ \ \ \ \ \ \ \ \ \  \ \ \ x\in \Omega
\end{cases}
\end{equation}
 \ }
where $0<T\leq +\infty$, $\Omega$ is an open and bounded domain of $\mathbb{R}^N$, $N\in\{1,2,3\}$ and $b_i\in H_0^2(\Omega):=H_0^1(\Omega)\cap H^2(\Omega)$. Here, $u_i(t)$ design the controls and  $y(t)$    the corresponding mild solution of the system (\ref{P3-sys obj}).

\medskip

In term of applications, equation like \eqref{P3-sys obj} may for instance describe the situation where some physical quantities (particles, energy,...) are transferred inside a system due to diffusion and convection processes, and the control can be seen as the velocity field that the quantity is moving with. For our special case, system \eqref{P3-sys obj} describes a deterministic Fokker–Planck equation for the time-dependent probability density ${ P\left(A,t\right)}$ of a stochastic variable $A$ of the Langevin equation, which enables us to study many types of fluctuations in physical and biological systems (see e.g.  \cite{krall}).

\medskip
 The goal of this paper is to study the problem of stability by an optimal control in the infinite time-horizon for the non-homogeneous bilinear system (\ref{P3-sys obj}). The main difficulty in solving a quadratic optimal control for general bilinear systems is the non-convexity of the cost function. In the case of a bounded control operator, the question of bilinear optimal control problem has been widely studied in the literature  (see \cite{Banks,Bradly,Alami,Li and Young,yahy,Boukhari R}). However, the modeling may give rise to the unboundedness aspect of the operator of control of the obtained bilinear model (see \cite{Addou,fokker,Berrahmoune 2009,clerin,Fle}), which is the case of equation \eqref{P3-sys obj} where the control is acting in the coefficient of the divergence term.  In \cite{Addou}, the authors  have considered the  homogeneous  version of the system \eqref{P3-sys obj} (i.e. $b_i=0,\; i=1,..,N$), for which they characterized the optimal control for a finite time horizon. Moreover, the author in  \cite{fokker,clerin} has studied the same problem in the presence of a time and state-dependent perturbation for finite time horizon.

\medskip

 The paper is organized as follows:  In Section \ref{sec2}  we give a preliminary.   In Section \ref{sec3}, 
 we first solve the optimal control problem in a finite time-horizon for the system \eqref{P3-sys obj}, and then proceed to the case of infinite time horizon in which context we give a  stabilization result by optimal control. Finally, in Section 4, we present a numerical example.
\section{Setting of the problem and some a priori estimates} \label{sec2}
Let us consider the following spaces: $ H=L^2(\Omega),\; V=H_0^1(\Omega), $ $V^*= H^{-1}(\Omega)$ and $U=L^2(0,T;\mathbb{R}^N)$, and let us introduce the following operators : 
		\begin{itemize}
\item {$A: V  \rightarrow V^*, \ y\mapsto \Delta y$ which is a linear continuous operator,}

\item {the linear continuous operator $B: V \rightarrow V^*$, is defined by $By=1_\Omega.\nabla y=\sum_{i=1}^{N}, \frac{\partial y}{\partial x_i}$, here $1_\Omega$ is the vector $ (1,\cdots,1)$.}
	\end{itemize}
 For all $u:=(u_i)_{1\leq i\leq N}\in U$ and $y\in V$, we have
\begin{equation*}
 	u.(By+b) :=u.(\nabla y+b) =\sum_{i=1}^{N} y_i \, \left(\frac{\partial y}{\partial x_i}+b_i \right)\cdot
		\end{equation*}
Thus the system (\ref{P3-sys obj}) can be rewritten in the form
\begin{equation}\label{P3-sys prin1}
\begin{cases}
\dot{y}(t)= A y(t) +u(t).(By(t)+b)  \\
    y(0)=y_0  \in H.
\end{cases}
\end{equation}	The quadratic cost function $J$ to be minimized is defined by
		\begin{equation}\label{P3 JT}
J(u)=  \int_0^{T} \Vert y(t)\Vert_H^2 dt + \frac{r}{2}\int_0^{T}\Vert u(t)\Vert_{\mathbb{R}^N}^2 dt,
 \end{equation} where $r>0$, $u\in U$ and $y(t)$ is the respective solution to system (\ref{P3-sys prin1}) .
 \medskip

 Then, the optimal control problem may be stated as follows
		
        \begin{equation}\label{P3prob opt}
    \left\{%
\begin{array}{ll}
\min J(u)  \\
    u\in U \\
\end{array}%
\right.
\end{equation}
 For the wellposedness of the system (\ref{P3-sys prin1}), let us consider the following system
\begin{equation}\label{P3 sys global}
\begin{cases}
\dot{y}(t)=Ay(t)+u(t).By(t)+f(t)\\
y(0)=y_0
\end{cases}
\end{equation}
where $f\in L^2(0,T;V^*)$, and let us introduce the following functional  space  $$W(0,T)=\{ \phi \in L^2(0,T,V)\ \ / \ \ \dot{\phi} \in L^2(0,T;V^*)\}\cdot$$
Now, we recall the following existence  result with some a priory estimates (see \cite{Addou,fokker,clerin}).

\begin{lem} \label{1}
For all $u\in U$, there exists a unique  solution $y$ of the system (\ref{P3 sys global}), which is such that $$y\in W(0,T) \cap L^{\infty}(0,T;H)\cdot$$
Moreover, the following estimates hold
    \begin{equation}\label{P3 lemma 1 1}
    \Vert y\Vert_{L^2(0,T;V)}\leq \frac{1}{\sqrt2}\Vert y_0\Vert_H+\Vert f\Vert_{L^2(0,T;V^*)}
    \end{equation}
    \begin{equation}\label{ P 3 estimate 7}
    \Vert y\Vert_{L^{\infty}(0,T;H)}\leq \Vert y_0\Vert_H+\sqrt2\Vert f\Vert_{L^2(0,T;V^*)}
    \end{equation}
     \begin{equation}\label{P3 ypoint}
     \begin{aligned}
    \Vert \dot{y}\Vert_{L^(0,T;V^*)}\leq  &\left(\frac{1}{\sqrt 2}\Vert y_0\Vert_H+\Vert f\Vert_{L^2(0,T;V^*)} \right) \left(\alpha +\sqrt 2 \Vert u\Vert_U \right)\\
    &\ +\Vert f\Vert_{L^2(0,T;V^*)}
    \end{aligned}
    \end{equation}
where $\alpha$ is such that $\Vert Az\Vert_{V^*}\leq \alpha\Vert z\Vert_V$, for all  $z\in V$.
\end{lem}
\section{Characterization of the optimal control} \label{sec3}

\subsection{Existence of an optimal control}

\begin{thm}\label{2}
For any $y_0\in H$,  the problem (\ref{P3prob opt}) has at least one solution.
\end{thm}
\begin{pf}
Since the set  $\{J(u)/ u\in U\}$ is not empty and is bounded from below,  it admits a lower bound $J^*$.\\
Let $ (u_n )_{n\in\mathbb{N}}$ be a minimizing sequence such that $J(u_n)\rightarrow J^*$.\\		
Then the sequence  $(u_n)$ is bounded, so it admits a sub-sequence  denoted by $(u_n)$ as well, which weakly converges to $u^*\in U$.\\
Let $(y_n)$  be the sequence of  solutions of (\ref{P3-sys prin1}) corresponding to $(u_n)$.
According to Lemma \ref{1}, the sequences $\begin{aligned}
		\Vert y_n(0)\Vert_H, \Vert y_n\Vert_{L^2(0,T;V)}, \ \Vert y_n\Vert_{L^{\infty}(0,T;V)}, \ \Vert \dot{y}_n\Vert_{L^2(0,T;V^*)},\end{aligned}
        $ $\begin{aligned}\Vert Ay_n\Vert_{L^2(0,T;V^*)}\end{aligned}$
				and
		$
		\Vert u_n.(By_n+b)\Vert_{L^2(0,T;V^*)}$ are bounded, so $(y_n)$ admits a sub-sequence, also denoted by $(y_n)$, such that
        $$y_n \rightharpoonup y^* \  weakly \ in \ L^2(0,T;V),$$
        $$y_n \rightharpoonup y^* \  weakly * \ \ in \ L^{\infty}(0,T;H),$$
        $$\dot{y}_n\rightharpoonup \dot{y}^* \ weakly \ in \ L^2(0,T;V^*).$$
        In addition to this, the linear operator 
        $$\mathbb{A}: L^2(0,T;V) \rightarrow L^2(0,T;V^*)$$$$ y \mapsto \mathbb{A} y$$
        is continuous, from which it follows that
        $$\mathbb{A}y_n\rightharpoonup \mathbb{A} y^*\ weakly\ in \ L^2(0,T;V^*)\cdot$$
        Then since the embedding $W(0,T)\rightarrow L^2(0,T;H)$ is compact, $(y_n)$  admits a sub-sequence, still denoted by $(y_n)$, for which we have
        \begin{equation}\label{P3 y_ny^*}
        y_n\rightarrow y^* \ strongly \ in \ L^2(0,T;H).
        \end{equation}
				Taking into account that the operator  $\mathbb{B}: L^2(0,T;H) \rightarrow L^2(0,T;V^*)$ is linear and continuous, we deduce that $$u_n.(\mathbb{B}y_n+b)\rightharpoonup u^*.(\mathbb{B}y^*+b) \ \ weakly \ in \ L^2(0,T;V^*).$$
				Now, by taking the  limit we deduce that
         \begin{equation*}
         \begin{cases}
         \dot{y^*}(t)=Ay^*(t)+u^*(t).(By^*(t)+b) \\
         y^*(0)=y_0
         \end{cases}
         \end{equation*}
         In other words, $y^*$ is the  solution of the system (\ref{P3-sys prin1}) corresponding to control  $u=u^*$.\\
Using that the norm  $\Vert \cdot \Vert_{L^2(H)}$ is lower semi-continuous, it follows from the strong convergence of the sequence $y_n$ to $y^*$ in $L^2(0,T;H)$  that

		\begin{equation}\label{P3 y inf Ja}
			\int_0^T\Vert y^*(t)\Vert_H^2dt\leq \liminf_{n \rightarrow +\infty} \int_0^T\Vert y_n(t)\Vert_H^2 dt.
		\end{equation}
		Since $R: u \mapsto \int_0^T\Vert u(t)\Vert_U^2 dt $ is convex and lower semi-continuous with respect to the weak topology, we have (see Corollary III.8 of \cite{Brezis})
		\begin{equation}\label{P3 u inf Ja}
			R(u^*)\leq \liminf_{n \rightarrow +\infty} R(u_n)\cdot
		\end{equation}
		Combining the  formulas (\ref{P3 y inf Ja}) and (\ref{P3 u inf Ja}) we deduce that
		\begin{equation*}
      	\begin{aligned}
   J(u^*) & =  \int_0^T\Vert y^*(t)\Vert_H^2dt +\frac{r}{2}\int_0^T\Vert u^*(t)\Vert_{\mathbb{R}^N}^2dt \\
            &\leq  \liminf_{n\rightarrow +\infty} \int_0^T\Vert y_n(t)\Vert_H^2dt+\frac{r}{2}\liminf_{n\rightarrow +\infty} \int_0^T\Vert u_n(t)\Vert_{\mathbb{R}^N}^2dt\\\\
            &\leq  \liminf_{n\rightarrow +\infty} J(u_n)\\\\
            &\leq J^*\cdot
            	\end{aligned}
		\end{equation*}
        We conclude that  $J(u^*)=J^*$, and so $u^*$ is a solution of the problem (\ref{P3prob opt}).
				
				\end{pf}
				
		\subsection{Expression of the optimal control for finite time-horizon}
 In this subsection, we will provide informations about the optimal control.
\begin{thm} \label{thm3}
For all $T>0$,  the problem (\ref{P3prob opt}) admits a solution  $u^*$ which is given by:
\begin{equation*}
u_i^*(t)=-\dfrac{1}{r}\langle \phi(t),\frac{\partial y^*(t)}{\partial x_i}+b_i\rangle_{V^*,V}, \, \forall \, i=1,...,N,
\end{equation*}
where $y^*$ is the solution of the system \eqref{P3-sys prin1} corresponding to $u^*$ and  $\phi$ is the  solution of the following adjoint system
\begin{equation}\label{P3 adjoin}
\left\{
\begin{array}{ll}
					\dot{\phi}(t)= -A \phi(t)+u^*(t). B\phi(t)-2y^*(t)\\
					\phi(T)=0
				\end{array}
\right.
			\end{equation}
\end{thm}

\begin{pf}
First, let us show that the mapping
$$
U \rightarrow C(0,T; H)$$
$$ u \mapsto y_u
$$
 is Fr\'echet differentiable and that its
derivative $z_h$ at $u\in U$, for a given $h \in U$, is the unique solution of the following system
\begin{equation}\label{P3 z_h}
\begin{cases}
\dot{z_h}(t)= Az_h(t)+u(t).Bz_h(t)+h(t).(By_u(t)+b)\\
z_h(0)=0
\end{cases}
\end{equation}
Let $u\in U$ and let $y_u$ be the corresponding solution of the system \eqref{P3-sys prin1}. We claim that the linear mapping $h \mapsto z_h$ is continuous. Indeed, using the  estimate (\ref{ P 3 estimate 7}) for the system (\ref{P3 z_h}), we can find  some $M>0$ such that
$$
\Vert z_h\Vert_{L^{\infty}(0,T;H)}\leq \sqrt2 M\Vert h\Vert_U.
$$
Let us denote by  $y_{u+h}$  the solution of the system (\ref{P3-sys prin1}) corresponding to $u+h$, and let $z_h$ be the solution of the system (\ref{P3 z_h}) corresponding to $h$. Taking $z=y_{h+u}-y_u-z_h$, we can see that $z$ is the solution of the following system

\begin{equation}\label{P3 z}
\begin{cases}
\dot{z}(t)= Az(t)+u(t).Bz(t)+h(t).B(y_{h+u}(t)-y_u(t))\\
z(0)=0
\end{cases}
\end{equation}
So, according to (\ref{ P 3 estimate 7}) in Lemma \ref{1}, the following estimates hold for some $K>0$
\begin{equation}\label{P3 estim 1}
\begin{aligned}
\Vert z\Vert_{L^{\infty}(0,T;H)} &\leq \sqrt{2}\Vert h.B(y_{h+u}(t)-y_u(t))\Vert_{L^2(0,T;V^*)}\\
&\leq K\Vert h\Vert_U\Vert y_{h+u}-y_u\Vert_{L^{\infty}(0,T;H)}.
\end{aligned}
\end{equation}
Let us set $w=y_{h+u}-y_u$. Then $w$ is the solution of the following system

\begin{equation}\label{P3 zw}
\begin{cases}
\dot{w}(t)= Aw(t)+u(t).Bw(t)+h(t).(By_{h+u}(t)+b)\\
w(0)=0
\end{cases}
\end{equation}
Applying Lemma \ref{1}, the following estimates hold for some $K_1, K_2>0$
\begin{equation}\label{P3 estim 2}
\begin{aligned}
\Vert w\Vert_{L^{\infty}(0,T;H)} &\leq \sqrt{2}\Vert h.(By_{h+u}+b)\Vert_{L^2(0,T;V^*)}\\
& \leq \Vert h\Vert_U \left(K_1\Vert y_{h+u}\Vert_{L^{\infty}(0,T;H)}+K_2 \right).
\end{aligned}
\end{equation}
Then using (\ref{P3 estim 1}) and (\ref{P3 estim 2}) and taking into account that the mapping $u \mapsto y_u$ is continuous, we conclude that for some $K_3>0$, we have
\begin{equation*}
\Vert z\Vert_{L^{\infty}(0,T;H)}\leq K_3\Vert h\Vert_U^2,
\end{equation*}
and hence the mapping $ u \mapsto y_u$ is Fr\'echet differentiable from $U$ to $ C(0,T; H)$, and that the derivative at $u\in U$ is given by the system (\ref{P3 z_h}).\\
Since the mappings $y\mapsto \Vert y\Vert_{L^2(0,T;H)}^2$ and $u\mapsto \Vert u\Vert_U^2$ are Fr\'echet differentiable, we  deduce that $u\mapsto J(u)$ is Fr\'echet differentiable as well, and we have
\begin{equation}\label{P3 J' sys}
\begin{cases}
D_uJ.h=\langle J'(u),h\rangle_U\\
D_uJ.h=\int_0^T\langle 2y(t),z_h(t) \rangle_H \, dt + r\int_0^T\langle u(t),h(t)\rangle_{\mathbb{R}^N} dt.
\end{cases}
\end{equation}
The well-posedness of the system (\ref{P3 adjoin}) is guaranteed by Lemma \ref{1}, after the following change of variables

\begin{equation}\label{P3 change of the variable}
\begin{cases}
q(t)=\phi(T-t)\\
g(t)=2y(T-t)\\
v(t)=u(T-t)\\
q(0)=\phi(T)=0.
\end{cases}
\end{equation}\\
Indeed, this leads to the following equivalent Cauchy problem:
\begin{equation*}
\begin{cases}
\dot{q}(t)= A q(t)-v(t). Bq(t)+g(t)\\
q(0)=0.
\end{cases}
\end{equation*}
Let $y$ and $\phi$   be the mild  solution of the systems (\ref{P3-sys prin1}) and (\ref{P3 adjoin}) respectively. Then we have
{\footnotesize
\begin{equation*}
\begin{aligned}
\int_0^T\langle  2y(t)&,z_h(t) \rangle_H dt  =\int_0^T\langle-\dot{\phi}(t)-A\phi(t)+u(t).B\phi(t),z_h(t)\rangle_{V^*,V} dt\\
& =-\int_0^T\langle\dot{\phi}(t),z_h(t)\rangle_{V^*,V}+\langle\phi(t),Az_h(t)+u(t).Bz_h(t)\rangle_{V^*,V} dt\\
 &=-\int_0^T\langle\dot{\phi}(t),z_h(t)\rangle_{V^*,V}+\langle\phi(t),\dot{z_h}(t)-h(t).(By(t)+b)\rangle_{V^*,V} dt\\
& =-\int_0^T\langle\dot{\phi_n}(t),z_h(t)\rangle_{V^*,V}+\langle\phi(t),\dot{z_h}(t)\rangle_{V^*,V}dt\\
&\ +\int_0^T\langle\phi(t),h(t).(By(t)+b)\rangle_{V^*,V}dt\\
& =-\left(\langle\phi(T),z_h(T)\rangle_{V^*,V}-\langle\phi(0),z_h(0)\rangle_{V^*,V}\right)\\
&\ +\int_0^T\langle\phi(t),h(t).(By(t)+b)\rangle_{V^*,V} dt.
\end{aligned}
\end{equation*} \ }
Since $\phi(T)=0$ and $z_h(0)=0$,  we conclude that
\begin{equation}\label{P3 3}
\begin{aligned}
\int_0^T\langle 2y(t),z_h(t)\rangle_H dt & = \int_0^T\langle \phi(t),h(t).(By(t)+b)\rangle_{V^*,V} dt\\
& =\int_0^T\langle(By(t)+b))^*\phi(t),h(t)\rangle_{\mathbb{R}^N} dt.
\end{aligned}
\end{equation}
Combining the formulas (\ref{P3 J' sys}) and (\ref{P3 3}) we deduce that

\begin{equation}\label{P1 J'}
\langle J'(u)(t), h(t)\rangle_{\mathbb{R}^N} =\langle (By(t)+b)^*\phi(t)+ru(t),h(t)\rangle_{\mathbb{R}^N}\cdot
\end{equation}
Hence  the solution of the problem (\ref{P3prob opt}) satisfies
\begin{equation*}
u_i^*(t)=-\dfrac{1}{r}\langle \phi(t),\frac{\partial y(t)}{\partial x_i}+b_i\rangle_{V^*,V},\ \  i=1,..N\cdot
\end{equation*}This achieve this proof.
\end{pf}
\subsection{Optimal control and strong stabilization}\ \\
Let us consider the following quadratic cost function $J$:
\begin{equation*}
J(u)=  \int_0^{+\infty} \Vert y(t)\Vert_H^2 dt + \frac{r}{2}\int_0^{+\infty}\Vert u(t)\Vert_{\mathbb{R}^N}^2 dt,
 \end{equation*} where $r>0, \ u\in U=L^2(0,+\infty; \mathbb{R}^N)$ and $y$ is the  corresponding mild solution of the system (\ref{P3-sys prin1}).\\
\medskip
The optimal control problem may be stated as follows
		
        \begin{equation}\label{P3prob opt infty}
    \left\{%
\begin{array}{ll}
\min J(u)  \\
    u\in U_{ad}=\{u \in L^2(0,+\infty;\mathbb{R}^N) \ \ / \  J(u)< +\infty \}
\end{array}
\right.
\end{equation}
Our goal in this part is to give a solution of the  problem \eqref{P3prob opt infty}. For this end, we consider the sequence of  controls $(u_n)$  solutions of the problem (\ref{P3prob opt}) on $[0,T_n]$ for  an increasing sequence  $T_n$ such that $T_n \rightarrow +\infty$. Let us denote by $y_n$  the  solution on $[0,T_n]$ of the system (\ref{P3-sys prin1}), and by  $\phi_n$ the solution of the adjoint system  (\ref{P3 adjoin}).

We have the following result.
\begin{thm} \label{thm4}
Let us consider the control $u^*=(u_i^*)$ defined by: 
		\begin{equation}\label{stab-opt}
		u_i^*(t)=-\dfrac{1}{r}\langle \phi(t),\frac{\partial y^*(t)}{\partial x_i}+b\rangle_{V^*,V} \ \ \ for\ \ i=1...N
		\end{equation}
		where $\phi$ is a weak limit value of $(\phi_n)$ in $L^2(0,+\infty;V)$ and $y^*$ is the corresponding solution of the system \ref{P3-sys prin1}. Then
\begin{itemize}
\item  $u^*$ is a solution of the problem (\ref{P3prob opt infty})
\item  $u^*$ guarantees the strong stabilization of  the system (\ref{P3-sys prin1}).
\end{itemize}

\end{thm}

\begin{pf}
		Let us first observe that $U_{ad}\ne\emptyset $, as here the solution of the system (\ref{P3-sys prin1}) corresponding to $u=0$ is exponentially stable.\\
		 Let $J_n$ be the functional (\ref{P3 JT}) in $[0,T_n],$ and let us define the following sequence of globally defined controls:
		\begin{equation*}
		v_n(t)=\left\{
		\begin{array}{ll}
		u_n(t), \ \ \ \hbox{if} \ \ \   t\leq T_n\\
		0, \ \ \  \hbox{if} \ \  \ t> T_n.
		\end{array}
		\right.
		\end{equation*}
Since $u_n\in L^2(0,T_n; \mathbb{R}^N), $ it follows that $v_n \in L^2(0,+\infty;\mathbb{R}^N)$.		
		Let us consider the mapping :  $$R: \ v\mapsto \frac{r}{2} \displaystyle \int_0^{+\infty}\Vert v(t)\Vert_{\mathbb{R}^N}^2dt.$$		
	Let $v\in U_{ad}$ be fixed. Since $ u_n $ is a solution of the problem (\ref{P3prob opt}) in $[0,T_n]$, it comes that		
		$$
R(v_n)=\frac{r}{2}\int_0^{T_n} \Vert u_n(t)\Vert^2_{\mathbb{R}^N}dt\leq J_n(u_n)\leq J_n(v)\leq J(v).
       $$
       Thus $R(v_n)$ is bounded and so is $v_n$. We deduce that the sequence $(v_n)$ admits  a subsequence, still denoted by $(v_n)$, which weakly converges   to $u^*\in L^p(0,+\infty)$.
			
			\medskip
			
Similarly to the proof of Theorem \ref{2}, we deduce that there exists a subsequence of $(v_n)$, (which  can be also denoted by $(v_n)$) such that
		\begin{equation*}
			y_{v_n}\rightarrow y^* \ strongly\  in \ L^2(0,+\infty;H),
		\end{equation*}
		where $y^*$ is the mild solution of the system (\ref{P3-sys prin1})  corresponding to $u^*$ in infinite time-horizon (i.e. $T=+\infty$).
       Then we conclude that
       \begin{equation}\label{P3 22}
       \lim_{n\rightarrow+\infty}\int_0^{T_n}\Vert y_{v_n}(t)\Vert_H^2dt=\int_0^{+\infty}\Vert y^*(t)\Vert_H^2 dt.
       \end{equation}
The continuity of the mapping $R$ implies the lower semi-continuity w.r.t to the weak topology (see Corollary III.8 in \cite{Brezis}). We deduce that
		\begin{equation}\label{P3 2222}
			R(u^*)\leq\liminf_{n \rightarrow +\infty} R(v_n) \cdot
            \end{equation}
		Observing that
		\begin{equation*}
			J_n(u_n)=\int_0^{+\infty}\Vert y_{v_n}1_{(0,T_n)}(t)\Vert_H^2dt+\frac{r}{2}\int_0^{+\infty}\Vert v_n(t)\Vert_{\mathbb{R}^N}dt,
		\end{equation*}
		we derive via (\ref{P3 22}) and (\ref{P3 2222})
		\begin{equation}\label{P1 J liminf}
			J(u^*)\leq \liminf_{n \rightarrow +\infty} (J_{n}(u_n)) \cdot
		\end{equation}
		Let us show that the sequence $(J_{n}(u_n))_{n\in \mathbb{N}}$ converges to $J(u^*)$. For this end, we will show that the sequence $(J_{n}(u_n))_{n\in \mathbb{N}}$ is increasing and upper bounded by $J(u^*)$. We have
\medskip
		$$J_{n}(u_n)\leq J_{n}(u_{n+1})\leq J_{{n+1}}(u_{n+1}) \  \ and\ \  J_{n}(u_n)\leq J_{n}(u^*)\leq J(u^*),
		$$
from which it comes
		\begin{equation}\label{P1 J sup}
			\lim_{n \rightarrow +\infty}J_{n}(u_n)\leq J(u^*)\cdot
		\end{equation}
		Combining (\ref{P1 J liminf}) and (\ref{P1 J sup}), we conclude that
$$
\lim_{n \rightarrow +\infty} J_{n}(u_n)= J(u^*)\cdot
$$
Keeping in mind that $u_n$  is the  solution of the problem (\ref{P3prob opt}) on $[0,T_n]$, we conclude that:
\begin{equation*}
\begin{aligned}
 J_n(u_n)-J(v) &=\int_0^{T_n} \left(\Vert  u_n(t)\Vert_{\mathbb{R}^N}^2+\Vert y_n(t)\Vert_H^2 \right) \, dt \\
 &\ \ -\int_0^{+\infty} \left(\Vert v(t)\Vert_{\mathbb{R}^N}^2+\Vert y_v(t)\Vert_H^2 \right)\, dt\\
&=J_n(u_n)-J_n(v)-\int_{T_n}^{+\infty}(\Vert v(t)\Vert_{\mathbb{R}^N}^2+\Vert y_v(t)\Vert_H^2)dt\\
&\leq 0\cdot
\end{aligned}
\end{equation*}
Thus letting $n\longrightarrow +\infty$, we get
		$$
J(u^*)-J(v)= \lim_{n \rightarrow +\infty}J_n(u_n)-J(v)\leq 0.
$$
This shows that $u ^ *$ is a solution of the problem (\ref{P3prob opt infty}).
Let  $\phi_n$  be the solution of the adjoint system (\ref{P3 adjoin})  corresponding to $u_n$. By the change of variables given by (\ref{P3 change of the variable}), $q_n$ is solution of the following system
\begin{equation*}
\begin{cases}
\dot{q}_n(t)=Aq_n(t)-v_n(t).Bq_n(t)+g_n(t)\\
q_n(0)=0

\end{cases}
\end{equation*}
So, by the estimate (\ref{P3 lemma 1 1}) in Lemma 1, we have
$\begin{aligned} \Vert \phi_n(T_n-.)\Vert_{L^2(0,T;V)}=\Vert q_n\Vert_{L^2(0,T;V)} \leq \Vert g_n\Vert_{L^2(0,T;H)}\end {aligned}$ and $\begin{aligned} \Vert g_n\Vert_{L^2(0,T;H)}=\Vert y_n(T_n-.)\Vert_{L^2(0,T;V)}\cdot\end {aligned}$ Then the boundedness of $y_n$ implies that of  $\phi_n$  in $L^2(0,+\infty; V)$. So, we can deduce that  the sequence $(\phi_n)$ admits  a subsequence, still denoted by $(\phi_n)$, which  weakly converges   to $\phi\in L^2(0,+\infty;V)$.\\
Using the fact that $u_n\rightharpoonup  u^*$   in $U$, $y_n\rightarrow y^*$  in $L^2(0,+\infty;V)$ and $\phi_n\rightharpoonup\phi$  in $L^2(0,+\infty;V)$,  we conclude by Theorem \ref{thm3} that 	
\begin{equation*}
		u_i^*(t)=-\dfrac{1}{r}\langle \phi(t),\frac{\partial y^*(t)}{\partial x_i}+b\rangle_{V^*,V}, \forall \,  i=1,..,N.
		\end{equation*}
		
Let us now show that this controls lead  to a strongly stable system in closed loop. For all $0< s< t<+\infty, $  we have
\begin{equation*}
\begin{aligned}
\vert \Vert y^*(t)\Vert^2-\Vert y^*(s)\Vert^2\vert &=\vert\int_s^t2\langle \dot{y}^*(r),y^*(r)\rangle_{V^*,V} dr\vert \\
&\leq  2\Vert \dot{y}^*\Vert_{L^2(s,t;V^*)}\left(\int_s^t\Vert y^*(r)\Vert_V^2dr \right)^{\frac{1}{2}}.
\end{aligned}
\end{equation*}
Then, according to estimate  (\ref{P3 ypoint}), we have for some $M_{y_0}> 0$
$$
\left| \Vert y^*(t)\Vert^2-\Vert y^*(s)\Vert^2\right|\leq M_{y_0} \left(\int_s^t\Vert y^*(r)\Vert_V^2dr \right)^{\frac{1}{2}}.
$$
Using the fact that $\int_0^{+\infty} \Vert u^*(t)\Vert^2dt<+\infty$, we deduce via (\ref{P3 lemma 1 1})  that $\int_0^{+\infty}\Vert y^*(t)\Vert^2_{V}<+\infty$.
Then we conclude that
$$
\Vert y^*(t)\Vert \rightarrow 0 \; \hbox{as} \; t\rightarrow +\infty.
$$

\end{pf}

\section{A numerical example}\
Here, we will present simulations in which we show numerically the strong stability of the optimal trajectory $y^*$ and we further compare numerically the optimal control w.r.t some controls $v$ in terms of energy consumption.\\
Let us consider the following parameters: $\Omega=(0,1)$, $b=5$, $y_0(x) =10x(1-x)$, $r=1$  and $T=8$.\\
Then reporting the states norm of the system for both controls $u^*$ and $v_1=0$ in Figure 1, we can see that $u^*$  performs slightly better than the zero control.
This tendency is confirmed in the table below regarding the states norm and the energy consumed by the system under the optimal control and constant controls. \\
\begin{figure}[ht]
\centering
\includegraphics[width=10cm, height=6cm]{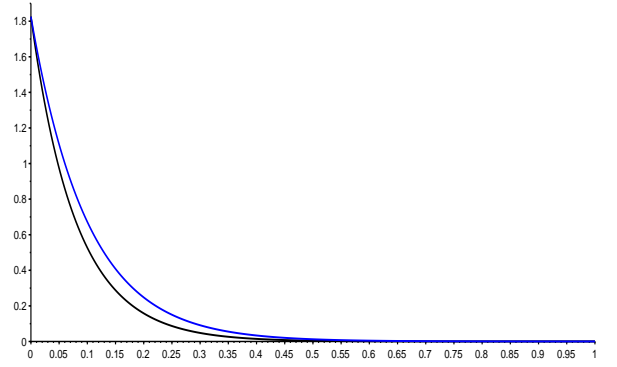}
\caption{Time-evolution of $\Vert y^*(t)\Vert$ under the optimal control $u^*$ (black line) and the zero control $v_1=0$ (blue line). }
\end{figure}

\begin{tabular}{c c c c }
\hline Time (t)  &  0.2 &  0.3  &0.6 \\
\hline  $\Vert y^*(t)\Vert$  &$15.8*10^{-2}$& $4.79*10^{-2}$ &$ 1.2*10^{-3}$\\
\hline $\Vert y_{v_1}(t)\Vert$  &$ 24.8*10^{-2} $&$ 9.12*10^{-2} $&$ 4.56*10^{-3} $ \\
\hline$ J(u^*)$ &  $14.68*10^{-2}$ & $14.77*10^{-2}$   & $14.77*10^{-2}$\\
\hline $J(v_1)$ & $16.38*10^{-2}$ & $16.63*10^{-2}$  & $16.67*10^{-2}$\\
\hline $J(v_2)$ & $4.19$ & $4.21$   & $4.22$\\
\hline
\end{tabular}

\begin{rem}
 Note that the stabilization problem of non-homogeneous distributed bilinear systems has been only considered for bounded control operator (see \cite{akouchi2003,bounit99,bounit2003,hamidi 17}). Thus the existing results from the above literature are not applicable as here, the control operator $B$ is unbounded. Moreover, even in the homogeneous case (i.e $b=0$), the existing results for unbounded operator $B$ (see \cite{ayadi,Berrahmoune 2010}) are not applicable as here, the operator $B$ is skew adjoint. In particular, the observation inequality is not verified. Now, if we formally consider the feedback control $v(t)$ used in \cite{ayadi,Berrahmoune 2010}, then we find $v(t)=0$  as $v(t)$ involves the term $\langle By,y\rangle$ which is null when $B$ is skew-adjoint.
\end{rem}
\section{Conclusion}
In this work, we studied the quadratic optimal control problem for a class of non-homogeneous bilinear Fokker-Planck equation. Both finite and infinite horizon cases are considered. It is further showed that the infinite horizon optimal control leads to a stabilized state of the system in closed-loop. This study provided a stabilization result which does not require the observation assumption. The result of Theorem 5 is promising. Indeed one can be inspired by it to investigate the optimal stabilization of a general unbounded bilinear system.

\end{document}